\documentclass[11pt]{article}
\usepackage{latex8}
\usepackage{times}
\usepackage{xspace}
\usepackage{paralist}
\usepackage{amsmath}

\usepackage{epsfig}

\newtheorem{theorem}{Theorem}[section]
\newtheorem{proposition}[theorem]{Proposition}
\newtheorem{lemma}[theorem]{Lemma}
\newtheorem{claim}[theorem]{Claim}
\def\squarebox#1{\hbox to #1{\hfill\vbox to #1{\vfill}}}
\newcommand{\qed}{\hspace*{\fill}
             \vbox{\hrule\hbox{\vrule\squarebox{.667em}\vrule}\hrule}\smallskip}
\newenvironment{proof}{\begin{trivlist}
\item[\hspace{\labelsep}{\bf\noindent Proof: }]
}{\qed\end{trivlist}}
\newcommand{\boldd}{{\bf d}}
\newcommand{\e}{{\rm e}}
\newcommand{\ex}{{\rm E}}
\newcommand{\ds}{{\rm d}s}
\newcommand{\dt}{{\rm d}t}
\newcommand{\eps}{\epsilon}
\newcommand{\smallspacing}{\baselineskip = 0.98\normalbaselineskip}
\newcommand{\normalspacing}{\baselineskip=.985\normalbaselineskip}
\newcommand{\goesto}{\rightarrow}
\def\PUR{{\rm PUR}}

\newcommand{\whp}{w.h.p.}

\newcommand{\remove}[1]{}

\begin{document}

\title{A Continuous-Discontinuous Second-Order Transition in the Satisfiability of Random 
Horn-SAT Formulas} 

\author{Cristopher Moore\thanks{
Department of Computer Science, University of New Mexico,
Albuquerque NM, USA, 
{\tt moore@cs.unm.edu}}
\and Gabriel Istrate\thanks{
CCS-DSS, Los Alamos National Laboratory, Los Alamos NM 87545, USA,
{\tt istrate@lanl.gov}}
\and Demetrios Demopoulos\thanks{
Archimedean Academy, 10870 SW 113 Place, Miami, FL, USA,
{\tt demetrios\_demopoulos@archimedean.org}}
\and Moshe Y. Vardi\thanks{
Department of Computer Science, Rice University,
Houston TX, USA,
{\tt vardi@rice.edu}}
}

\maketitle
\thispagestyle{empty}

\begin{abstract} 
We compute the probability of satisfiability of a class of 
random Horn-SAT formulae, motivated by a connection 
with the nonemptiness problem of finite tree automata.
In particular, when the maximum clause length is 3, 
this model displays a curve in its parameter space along which 
the probability of satisfiability is discontinuous, 
ending in a second-order phase transition where it becomes continuous.  
This is the first case in which a phase transition of this type has been 
rigorously established for a random constraint satisfaction problem. 
\end{abstract}

\section{Introduction} 

In the past decade, phase transitions, or {\em sharp thresholds}, have
been studied intensively in combinatorial problems.  Although the idea 
of thresholds in a combinatorial context was introduced as early as 1960~\cite{ER60}, 
in recent years it has been a major subject
of research in the communities of theoretical computer science,
artificial intelligence and statistical physics.  Phase transitions 
have been observed in numerous combinatorial problems in which they 
the probability of satisfiability jumps from $1$ to $0$ when the density of constraints 
exceeds a critical threshold. 

The problem at the center of this research is, of course, 3-SAT.
An instance of 3-SAT is a Boolean formula, consisting of a conjunction
of clauses, where each clause is a disjunction of three literals. The goal 
is to find a truth assignment that satisfies all the clauses and thus the entire formula. 
The {\em density} of a 3-SAT instance is the ratio of the number of clauses to the number of
variables.  We call the number of variables the {\em size\/} of the instance.  
Experimental studies~\cite{CA96,SML96,SK96}
show a dramatic shift in the probability of satisfiability 
of random 3-SAT instances, from 1 to 0, located at a critical density $r_c \approx 4.26$.
However, in spite of compelling arguments from statistical physics~\cite{MZ02,MPZ02}, 
and rigorous upper and lower bounds on the threshold if it exists~\cite{DBM00,HajiaSorkin,KKL}, 
there is still no mathematical proof that a phase transition takes place at that density.
\remove{
A limitation on all the experimental studies is
imposed by the inherent difficulty of the problem, especially around the
phase-transition point. We can only study instances of limited size 
(usually up to a few hundred) before the problems get too hard to be solved 
in reasonable time using available computational resources.
}
For a view variants of SAT the existence and location of 
phase transitions have been established rigorously, in particular for
2-SAT, 3-XORSAT, and 1-in-$k$ SAT~\cite{ACIM01,DM02,CDMM03,CR92,Goe96}.  

In this paper we prove the existence of a more elaborate type of 
phase transition, where a curve of discontinuities in a two-dimensional parameter 
space ends in a {\em second-order} transition where the probability of satisfiability 
becomes continuous.  We focus on a particular variant of 3-SAT, namely 
Horn-SAT.  A Horn clause is a disjunction of literals of which {\em at most one} is positive, 
and a Horn-SAT formula is a conjunction of Horn clauses.  Unlike 3-SAT, Horn-SAT is a 
tractable problem; the complexity of the Horn-SAT is linear in
the size of the formula~\cite{DG84}.  This tractability allows
one to study random Horn-SAT formulae for much larger input sizes 
that we can achieve using complete algorithms for 3-SAT.

An additional motivation for studying random Horn-SAT comes from the fact
that Horn formulae are connected to several other areas of computer science
and mathematics~\cite{Mak87}.  In particular, Horn formulae are connected to 
automata theory, as the transition relation, the starting state, and the set 
of final states of an automaton can be described using Horn clauses.
For example, if we consider automata on binary trees, 
then Horn clauses of length three can be used to describe its 
transition relation, while Horn clauses of length one can describe
the starting state and the set of the final states of the automaton.  
(we elaborate on this below).  
Then the question of the emptiness of
the language of the automaton can be translated to a question about the 
satisfiability of the formula.  Since automata-theoretic techniques have
recently been applied in automated reasoning~\cite{VW86a,VW86b},
the behavior of random Horn formulae might shed light on these applications.

Threshold properties of random Horn-SAT problems have recently been  
actively studied. The probability of satisfiability of random 
Horn formulae in a \emph{variable}-clause-length 
model was fully characterized in~\cite{Ist02,Ist04},
where it was shown that random Horn formulae have a {\em coarse} rather 
than a sharp satisfiability threshold, meaning that the
problem does not have a phase transition in this model. 
The variable-clause-length model used there is ideally suited to studying Horn
formulae in connection with knowledge-based systems~\cite{Mak87}.
Bench-Capon and Dunne~\cite{BCD01} studied
a \emph{fixed}-clause-length model, in which each Horn clause has 
precisely $k$~literals, and proved a sharp threshold with respect to 
assignments that have at least $k-1$ variables assigned to be true.

\begin{figure}[tbh]
\begin{center}
\begin{minipage}[t]{2.5 in}
\begin{center}
\psfig{figure=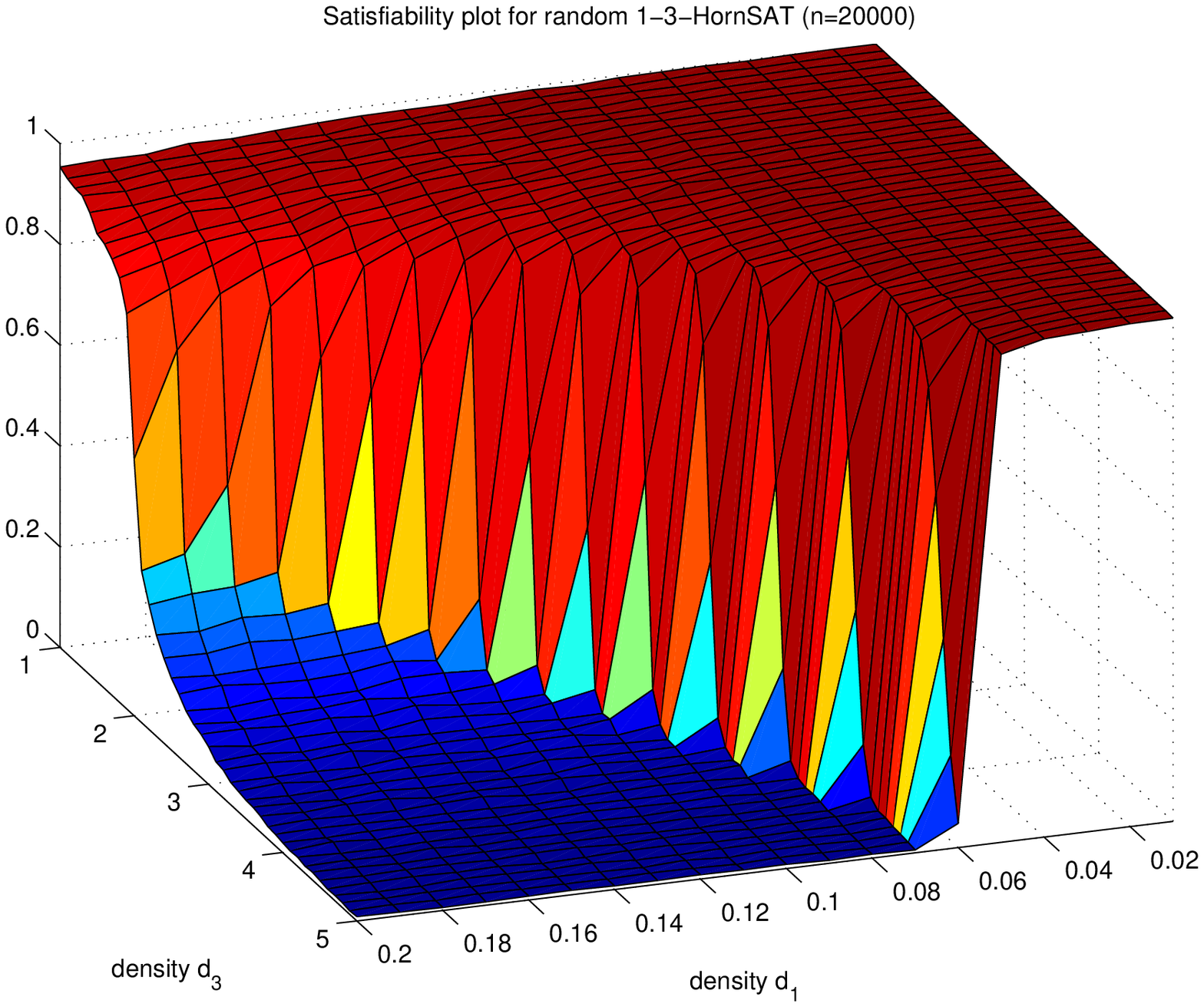,height= 70 truemm,width=65 truemm}
\end{center}
\end{minipage}
\begin{minipage}[t]{2.8 in}
\begin{flushright}
\hspace{0.5in}\psfig{figure=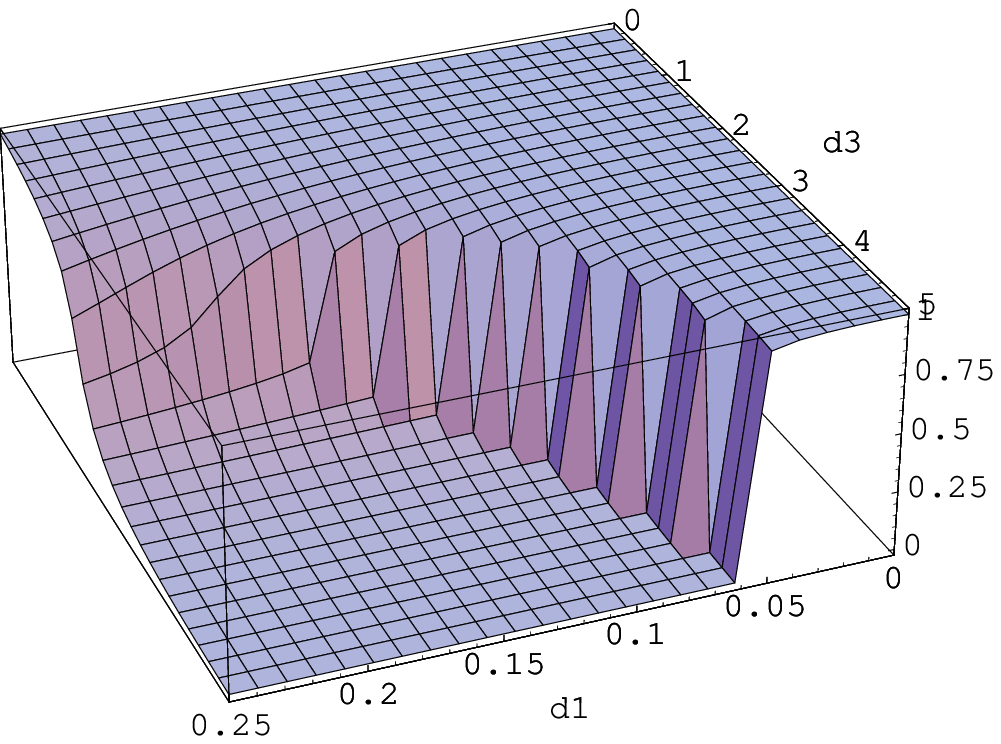,height= 70 truemm,width= 65 truemm} 
\end{flushright}
\end{minipage}
\caption{Satisfiability probability of a random 1-3-Horn
formula of size 20000.  Left, the experimental results of~\cite{DV05}; 
right, our analytic results.}
\label{fig:1-3}
\end{center}
\end{figure}

Motivated by the connection between the automata emptiness problem and
Horn satisfiability, Demopoulos and Vardi~\cite{DV05} studied 
the satisfiability of two types
of fixed-clause-length random Horn-SAT formulae.  They considered 
1-2-Horn-SAT, where formulae consist of clauses of length one and two 
only, and 1-3-Horn-SAT, where formulae consist of clauses of length one 
and three only.  These two classes can be viewed as the Horn analogue of
2-SAT and 3-SAT. For 1-2-Horn-SAT, they showed experimentally 
that there is a coarse transition
(see Figure~\ref{1-2-0.1}), which can be explained and analyzed
in terms of random digraph connectivity~\cite{Kar90}. 
The situation for of 1-3-Horn-SAT is less clear cut.  On one hand, 
recent results on random undirected hypergraphs~\cite{DN05} fit the 
experimental data of~\cite{DV05} quite well.  On the other, a scaling analysis of the data 
suggested that transition between the mostly-satisfiable and mostly-unsatisfiable regions
(the ``waterfall'' in Figure~\ref{fig:1-3}) is steep but continuous, rather than a step function.  
It was therefore not clear if the model exhibits a phase transition, in spite of 
experimental data for instances with tens of thousands of variables.

In this paper we generalize the fixed-clause-length model of
\cite{DV05} and offer a complete analysis of the probability of 
satisfiability in this model.  For a finite $k>0$ and 
a vector $\boldd$ of~$k$ nonnegative real numbers 
${d_1,d_2,\ldots,d_k}$, $d_1 < 1$, let the random Horn-SAT formula 
$H^k_{n,\boldd}$ be the conjunction of
\begin{itemize}
\item a single negative literal $\overline{x}_1$,
\item $d_1 n$ positive literals chosen uniformly without replacement from $x_2,\ldots,x_n$, and
\item for each $2 \le j \le k$, $d_j n$ clauses chosen uniformly from the 
$j{n \choose j}$ possible Horn clauses with $j$ variables where one literal is positive.
\end{itemize}
Thus, the classes studied in~\cite{DV05}
are $H^2_{n,d_1,d_2}$ and $H^3_{n,d_1,0,d_3}$ respectively.

With this model in hand, we settle
the question of sharp thresholds for 1-3-Horn-SAT.  In particular, we show that
there are sharp thresholds in some regions of the $(d_1,d_3)$ plane in the 
probability of satisfiability, although not from $1$ to $0$.  We start with the following 
general result for the $H^k_{n,\boldd}$ model.

\begin{theorem}\label{thm:main}  
Let $t_0$ be the smallest positive root of the equation
\begin{equation}
\label{eq:1}
 \ln \frac{1-t}{1-d_1} + \sum_{j=2}^k d_j t^{j-1} = 0 \enspace . 
\end{equation}
Call $t_0$ {\em simple} if it is not a double root of~\eqref{eq:1}, i.e., if the derivative of the left-hand-side of~\eqref{eq:1} with respect to $t$ is nonzero at $t_0$.  If $t_0$ is simple, 
the probability that a random formula from $H^k_{n,\boldd}$ is satisfiable in the limit 
$n \to \infty$ is 
\begin{equation}
\label{eq:prsat}  
\Phi(\boldd):= 
\lim_{n\goesto \infty} \Pr[H^k_{n,\boldd} \mbox{ is satisfiable}] = 
\frac{1-t_0}{1-d_1} 
\enspace . 
\end{equation} 
\end{theorem} 

Specializing this result to the case $k=2$ yields an exact formula 
that matches the experimental results in~\cite{DV05}: 
\begin{proposition}\label{prop:2}
The probability that a random formula from 
$H^2_{n,d_1,d_2}$ is satisfiable in the limit $n \to \infty$ is 
\begin{equation}
\label{eq:prop2}
\Phi(d_1,d_2):=
\lim_{n\goesto \infty}  \Pr[H^2_{n,d_1,d_2} \mbox{ is satisfiable}] = 
-\frac{W\bigl(-(1-d_1)d_2e^{-d_2}\bigr )}{(1-d_1)d_2} 
\enspace . 
\end{equation}
\end{proposition} 
Here $W(\cdot)$ is Lambert's function, defined as the principal root of 
the equation $W(x)e^{W(x)}=x$.

For the case $k=3$ and $d_2=0$, we do not have a closed-form
expression for the probability of satisfiability, though numerically Figure~\ref{fig:1-3} shows
a very good fit to the experimental results of~\cite{DV05}. 
In addition, we find an interesting phase transition behavior 
in the $(d_1,d_3)$ plane, described by the following proposition.
\begin{proposition}\label{prop:3}
The probability of satisfiability $\Phi(d_1,d_3)$ that a random formula from $H^3_{n,d_1,0,d_3}$ 
is satisfiable is continuous for $d_3 < 2$ and discontinuous for $d_3 > 2$.  Its discontinuities are given by a curve $\Gamma$ in the $(d_1,d_3)$ plane described by the equation 
\begin{equation}
\label{eq:gamma}
d_1 
 = 1 - \frac{\exp\left( \frac{1}{4} \left( \sqrt{d_3} - \sqrt{d_3-2} \right)^2 \right)}{d_3 - \sqrt{d_3 (d_3-2)}} 
 \enspace .
\end{equation}
This curve consists of the points $(d_1,d_3)$ at which $t_0$ is a double root of~\eqref{eq:1}, and ends at the critical point
\[ (1-\sqrt{\e}/2,2) = (0.1756...,2) \enspace . \]
\end{proposition} 

\begin{figure*}[tbh]
\begin{minipage}[t]{3.2  in}
\begin{center}
\psfig{figure=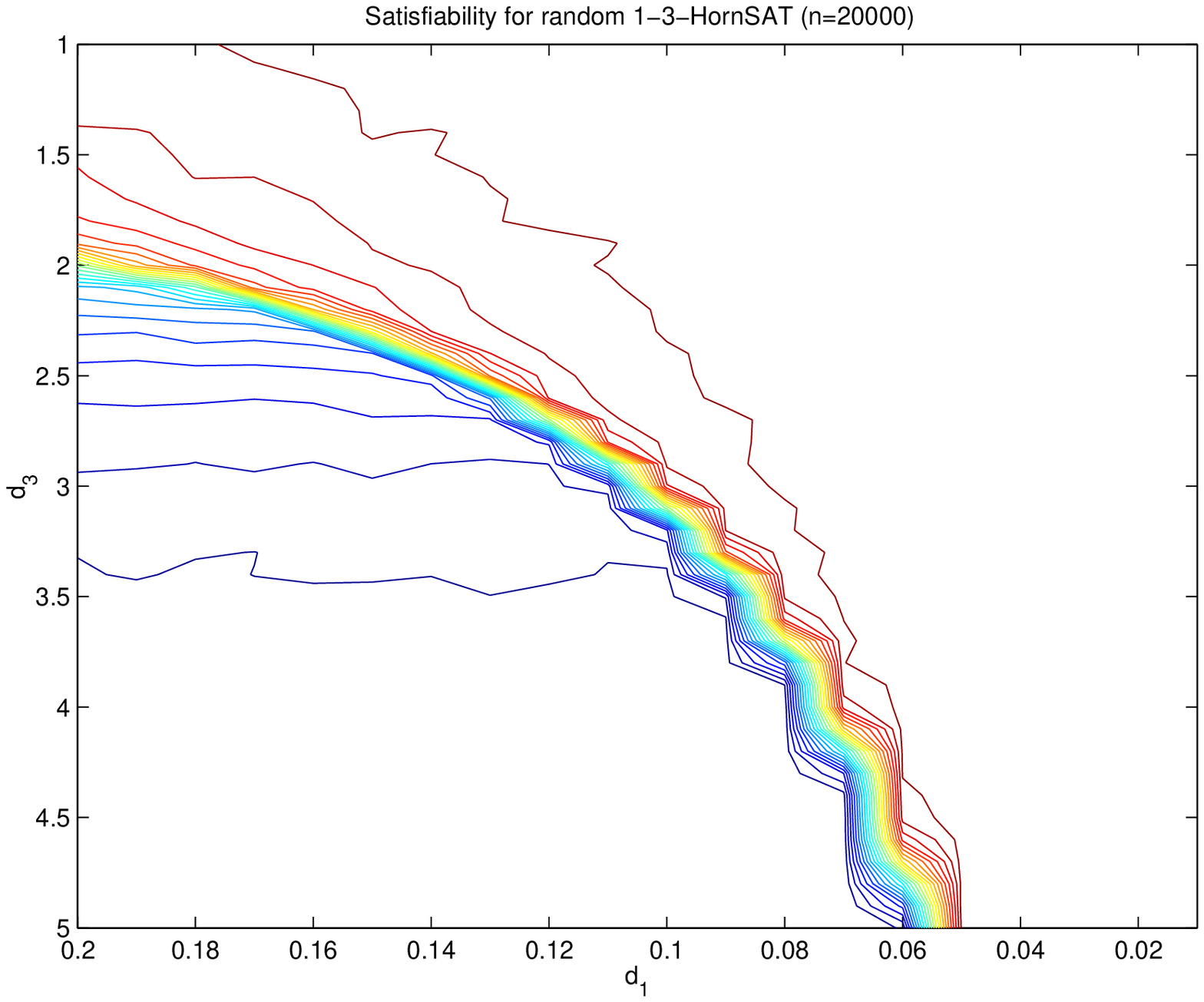, height= 45 truemm, width= 50 truemm}
\end{center}
\end{minipage}
\begin{minipage}[t]{3.3 in}
\begin{flushright}
\hspace{0.5in}\psfig{figure=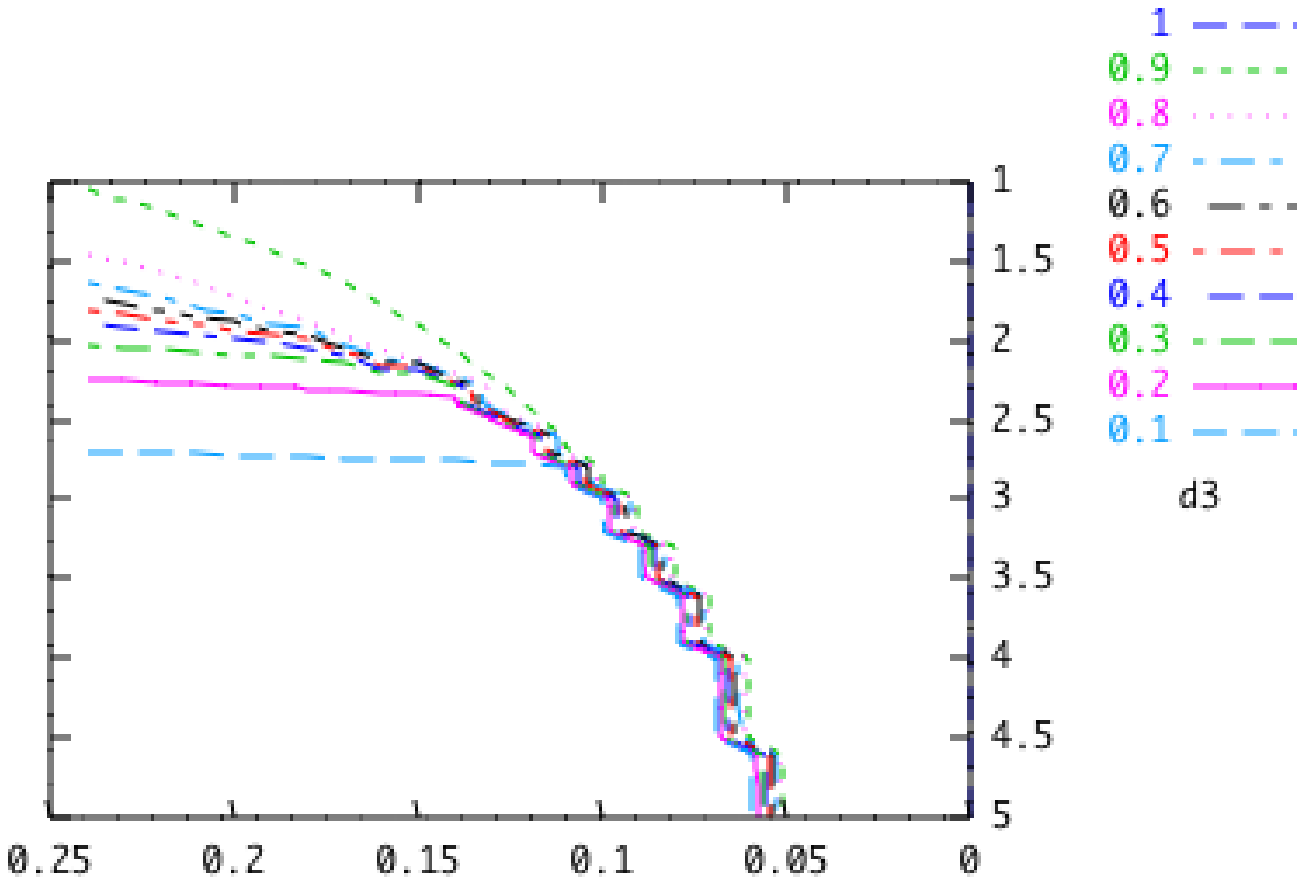,height= 58 truemm, width= 60 truemm}
\end{flushright}
\end{minipage}
\caption{Contour plots.  Left, the experimental results of~\cite{DV05}.  Right, our analytical results.}
\label{fig:contour}
\end{figure*}

The curve $\Gamma$ of discontinuities described in Proposition~\ref{prop:3} can be seen in the right part of Figure~\ref{fig:1-3}.  The drop at the ``waterfall'' decreases as we approach the critical point $(0.1756...,2)$, where the probability of satisfiability becomes continuous (although its derivatives at that point are infinite).  We can also see this in Figure~\ref{fig:contour}, which shows a contour plot; the contours are bunched at the discontinuity, and ``fan out'' at the critical point.  In both cases our calculates closely match the experimental results of~\cite{DV05}.

\remove{
Unlike much of the work on phase transitions in computer science,
rather than varying a single parameter, we explore a ``phase plane'' of
two parameters, namely the density of 3-clauses and unit clauses in the
formula.  In this plane we see a very interesting behavior.  There is a
curve along which $\Pr[\mbox{SAT}]$ jumps
discontinuously from one value to another; however, this curve ends at
a particular point, at which $\Pr[\mbox{SAT}]$ becomes continuous
(although its derivatives at that point are infinite).  
The endpoint $d_1=(1-\sqrt{\e}/2,2)$, $d_3=2$ of the discontinuity 
curve marks a transition between two 
regimes: for $d_3<2$ the probability of satisfiability is continuous, while
for $d_3>2$  the probability of satisfiability becomes discontinuous.
}

In statistical physics, we would say that $\Gamma$ is a curve of 
{\em first-order} transitions, in which the order parameter is 
discontinuous, ending in a {\em second-order} transition at the tip of 
the curve, at which the order parameter is continuous, but has a 
discontinuous derivative (see e.g.~\cite{BDFN92}).  
A similar transition takes place in the Ising
model, where the two parameters are the temperature $T$ and the
external field $H$; the magnetization is discontinuous at the
line $H=0$ for $T < T_c$ where $T_c$ is the transition temperature,
but there is a second-order transition at $(T_c,0)$ and the
magnetization is continuous for $T \ge T_c$.

To our knowledge, this is the first time that a
continuous-discontinuous transition of this type has been
established rigorously in a model of random constraint satisfaction
problems.  We note that a similar phenomenon is believed to take place 
for $(2+p)$-SAT model at the triple point $p=2/5$;
here the order parameter is the size of the backbone, i.e., 
the number of variables that take fixed 
values in every truth assignment~\cite{AKKK01,MZKST99}.  Indeed, 
in our model the probability of satisfiability is closely related to the 
size of the backbone, as we see below.

\remove{
Theorem~\ref{thm:main} allows us to settle the apparent 
contradiction  in~\cite{DV05}
between the conclusion sugegsted by scaling anaysis and those 
suggested by the random undirected hypergraphs approximation.
As the scaling analysis suggests, there is a sharp theshold for
satisfiability, but as our analysis here shows, the transition 
is not from 1 to 0.
}

\section{Horn-SAT and Automata}

Our main motivation for studying the satisfiability of Horn
formulae is the unusually rich type of phase transition described above, 
and the fact that its tractability allows us to perform experiments on 
formulae of very large size.  However, the original motivation of~\cite{DV05} that led 
to the present work is the fact that Horn formulae can be 
used to describe finite automata on words and trees.  

A finite automaton $A$ is a 5-tuple 
$A=(S,\Sigma,\delta,s,F)$, where $S$ is a finite set of states, $\Sigma$
is an alphabet, $s$ is a starting state, $F\subseteq S$ is the set of 
final (accepting) states and $\delta$ is a transition relation. 
In a word automaton, $\delta$ is a function from $S\times\Sigma$ to
$2^S$, while in a binary-tree automaton $\delta$ is a function from
$S\times\Sigma$ to $2^{S\times S}$. 
Intuitively, for word automata $\delta$ provides a set of
successor states, while for binary-tree automata $\delta$
provides a set of successor state pairs.
A run of an automaton on a word
$a=a_1a_2\cdots a_n$ is a sequence of states $s_0s_1\cdots s_n$ such
that $s_0=s$ and $(s_{i-1},a_i,s_i)\in\delta$. 
A run is succesful if
$s_n\in F$: in this case we say that $A$ accepts the word $a$. A run
of an automaton on a binary tree $t$ labeled with letters from
$\Sigma$ is a binary tree $r$ labeled with states from $S$ such that
root$(r)=s$ and for a node $i$ of $t$,
$(r(i),t(i),r(\mbox{left-child-of-}i),r(\mbox{right-child-of-}i))\in\delta$.
Thus, each pair in $\delta(r(i),t(i))$ is a possible labeling of
the children of~$i$.
A run is succesful if for all leaves $l$ of $r$, $r(l)\in F$: in this
case we say that $A$ accepts the tree $t$. The language $L(A)$ of a
word automaton $A$ is the set of all words $a$
for which there is a successful run of $A$ on $a$.  Likewise, the
language $L(A)$ of a tree automaton $A$ is the set of all trees $t$
for which there is a successful run of $A$ on $t$.
An important question in automata theory, which is also of
great practical importance in the field of formal verification
\cite{VW86a,VW86b}, is: given an automaton $A$, is $L(A)$ non-empty? 
We now show how the problem of non-emptiness of automata 
languages translates to Horn satisfiability. Thus, getting a
better understanding of the satisfiability of Horn formulae
would tell us about the expected answer to automata nonemptiness
problems.

Consider first a word automaton
$A=(S,\Sigma,\delta,s_0,F)$. Construct a Horn formula $\phi_A$ over the 
set $S$ of variables as follows: 
create a clause $(\bar{s_0})$,
for each $s_i\in F$ create a clause $(s_i)$,
for each element $(s_i,a,s_j)$ of $\delta$ create a clause
$(\bar{s_j}, s_i)$,
where $(s_i,\cdots,s_k)$ represents the clause $s_i\vee\cdots\vee s_k$ 
and $\bar{s_j}$ is the negation of $s_j$.
Similarly to the word automata case, we can show how to construct a Horn
formula from a binary-tree automaton. Let $A=(S,\Sigma,\delta,s_0,F)$ be
a binary-tree automaton. Then we can construct a Horn formula $\phi_A$ 
using the construction above with the only difference that since 
$\delta$ in this case is a function from $S\times\{\alpha\}$ to 
$S\times S$, for each element $(s_i,\alpha,s_j,s_k)$ of $\delta$ we 
create a clause $(\bar{s_j},\bar{s_k},s_i)$. 

\begin{proposition}
{\rm~\cite{DV05}}
Let $A$ be a word or binary tree automaton and $\phi_A$ the Horn formula
constructed as described above. Then $L(A)$ is non-empty if and only if 
$\phi_A$ is unsatisfiable.
\end{proposition}

\section{Main Result} 
\label{main}

\begin{figure}[t]
\small{
\begin{center}
\noindent \fbox{
\begin{minipage}[p]{3.75in}
\small{
\smallspacing
\medskip
\noindent
{\bf Algorithm $\PUR$:}
\begin{asparaenum}
\item while ($\phi$ contains positive unit clauses) 
\item \hspace{5mm} choose a random positive unit clause $x$
\item \hspace{5mm} remove all  other clauses in which $x$ occurs positively
in $\phi$  
\item \hspace{5mm} shorten all clauses in which $x$ appears negatively
\item \hspace{5mm} label $x$ as ``implied'' and call the algorithm recursively.
\item if no contradictory clause was created 
\item \hspace{5mm} accept $\phi$  
\item else 
\item \hspace{5mm} reject $\phi$. 
\end{asparaenum}
}
\end{minipage}
}
\end{center}
\caption{\textbf{Positive Unit Resolution.}}
\label{pur}
}
\end{figure}

\normalsize
\normalspacing

Consider the positive unit resolution algorithm $\PUR$ applied to a random formula $\phi$ (Figure~\ref{pur}). The proof of Theorem~\ref{thm:main} follows immediately from the following theorem, which establishes the size of the backbone of the formula with the single negative literal $\overline{x}_1$ removed: that is, the set of positive literals implied by the positive unit clauses and the clauses of length greater than $1$.  Then $\phi$ is satisfiable as long as $x_1$ is not in this backbone.

\begin{lemma}
\label{lem:gen}
Let $\phi$ be a random Horn-SAT formula $H^k_{n,\boldd}$ with $d_1 > 0$.  Denote by $t_0$ the smallest positive root of Equation~\eqref{eq:1}, and suppose that $t_0$ is simple.  Then, for any $\eps > 0$, the number $N_{\boldd,n}$ of implied positive literals, including the $d_1 n$ initially positive literals, satisfies \whp\ the inequality 
\begin{equation}
\label{eq:implied} 
(t_0-\eps) \cdot n < N_{\boldd,n} < (t_0+\eps) \cdot n, 
\end{equation}  
\end{lemma}

\begin{proof}
First, we give a heuristic argument, analogous to the branching process argument for the size of the giant component in a random graph.  
The number $m$ of clauses of length $j$ with a given literal $x$ as their implicate (i.e., in which $x$ appears positively) is Poisson-distributed with mean $d_j$.  If any of these clauses have the property that all their literals whose negations appear are implied, then $x$ is implied as well.  In addition, $x$ is implied if it is one of the $d_1 n$ initially positive literals.  Therefore, the probability that $x$ is {\em not} implied is the probability that it is not one of the initially positive literals, and that, for all $j$, for all $m$ clauses $c$ of length $j$ with $x$ as their implicate, at least one of the $j-1$ literals whose negations appear in $c$ is not implied.  Assuming all these events are independent, if $t$ is the fraction of literals that are implied, we have
\begin{align*}
 1-t & = (1-d_1) \prod_{j=2}^k 
\left( \sum_{m=0}^{\infty} \frac{\e^{-d_j} d_j^{m}}{m!} (1-t^{j-1})^{m} \right) \\ 
 & = (1-d_1) \prod_{j=2}^k \exp(-d_j t^{j-1})) 
 = (1-d_1) \,\exp\left( - \sum_{j=2}^k d_j t^{j-1}\right) 
\end{align*}
yielding Equation~\eqref{eq:1}.

To make this rigorous, we use a standard technique for proving results about 
threshold properties: analysis of algorithms via differential equations~\cite{Wor95} (see~\cite{achreview} for a review).  We analyze the while loop of $\PUR$ shown in Figure~\ref{pur}; specifically, we view $\PUR$ as working in {\em stages}, indexed by the 
number of literals that are labeled ``implied.'' 
After $T$ steps of this process, $T$ variables are labeled as implied. 
At each stage the resulting formula consists of a set of 
Horn clauses of length $j$ for $1 \le j \le k$ on the $n-T$ unlabeled variables.  
Let the number of distinct clauses of length $j$ in this formula be $S_j(T)$; 
we rely on the fact that, at each stage $T$, conditioned on the values of $S_j(T)$ 
the formula is uniformly random.  This follows from an easy induction argument
which is standard for problems of this type (see e.g.~\cite{Ist04}).

\remove{
The technical result on which our analysis rests is the following 
``uniformity lemma'': 
\begin{claim}
Conditional on the values of $S_j(T)$,
the formula at stage $T$ is random and uniform.   
\end{claim} 
}

Now, the variables appearing in the clauses present at stage $T$ are chosen 
uniformly from the $n-T$ remaining variables, so the probability that the chosen variable $x$ 
appears in a given clause of length $j$ is $j/(n-T)$, and the probability 
that a given clause of length $j+1$ is shortened to one of length $j$ 
(as opposed to removed) is $j/(n-T)$.  A newly shortened clause is 
distinct from all the others with probability $1-o(1)$ unless $j=1$, 
in which case it is distinct with probability $(n-T-S_1)/(n-T)$.  
Finally, each stage labels the variable in one of the $S_1(T)$ unit clauses as implied.  
Thus the expected effect of each step is
\begin{eqnarray*}
\ex[S_j(T+1)] & =& S_j(T) + j \frac{S_{j+1}(T) - S_j(T)}{n-T} + o(1)
\quad \mbox{for all } 2 \le j \le k \\
\ex[S_1(T+1)] & =& S_1(T) + \left( \frac{n-T-S_1}{n-T} \right) \left( \frac{S_2(T)}{n-T} \right) - 1 + o(1)
\end{eqnarray*}

\remove{
\begin{eqnarray*}
S_j(T+1)=S_j(T) & + & \Delta_{j+1}(T) - \Xi_{j}(T), 
\quad \mbox{for all } 2 \le j \le k \\
S_1(T+1)=S_1(T) & - & 1+\Delta_{2}(T) - \Xi_{1}(T),
\end{eqnarray*}
where $S_{k+1}(T) = 0$ and parameters $\Delta_j(T+1), \Xi_j(T+1)$ have the following distributions: 
\begin{eqnarray*}
\Delta_j(T)& \stackrel{D}{=} & B(S_j(T),(j-1)/(n-T)), \\
\Xi_j(T)& \stackrel{D}{=} & B(S_j(T),j/(n-T)).
\end{eqnarray*}
for $j\geq 2$, and 
\begin{eqnarray*}
\Xi_1(T+1)\stackrel{D}{=} B(S_1(T),(n-T-S_1(T))/(n-T)).
\end{eqnarray*}
}

\noindent 
Setting $T=t \cdot n$ and $S_j(T)=s_j(t) \cdot n$, we rescale this to form a system of differential equations:
\begin{eqnarray}
\frac{\ds_j}{\dt} & = & j \,\frac{s_{j+1}(t)-s_j(t)}{1-t}
\quad \mbox{for all } 2 \le j \le k \nonumber \\
\frac{\ds_1}{\dt} & = & \left( \frac{1-t-s_1(t)}{1-t} \right) \left( \frac{s_2(t)}{1-t} \right) - 1
\enspace . 
\label{eq:sys}
\end{eqnarray}
Then Wormald's Theorem tells us that, for any constant $\delta > 0$, for all $t$ such that $s_1 > \delta$, \whp\ we have $S_j(t \cdot n) = s_j(t) \cdot n + o(n)$ where $s_j(t)$ is the solution to the system~\eqref{eq:sys}.  With the initial conditions $s_j(0) = d_j$ for $1 \le j \le k$, a little work shows that this solution is
\begin{eqnarray}
s_j(t) & = & (1-t)^j \,\sum_{\ell = j}^k {\ell-1 \choose j-1} \,d_\ell t^{\ell-j}
\quad \mbox{for all } 2 \le j \le k \nonumber \\
s_1(t) & = & 1-t - (1-d_1) \exp\left(-\sum_{j=2}^k d_j t^{j-1} \right) \enspace . 
\label{eq:soln}
\end{eqnarray}

Note that $s_1(t)$ is continuous, $s_1(0)=d_1>0$, and $s_1(1) < 0$ since $d_1 < 1$.  Thus $s_1(t)$ has at least one root in the interval $(0,1)$.  Since $\PUR$ halts when there are no unit clauses, i.e., when $S_1(T)=0$, we expect the stage at which it halts.  Thus the number of implied positive literals, to be $T=t_0 n + o(n)$ where $t_0$ is the smallest positive root of $s_1(t)=0$, or equivalently, dividing by $1-d_1$ and taking the logarithm, of Equation~\eqref{eq:1}.

However, the conditions of Wormald's theorem do not allow us to run the differential equations all the way up to stage $t_0 n$. We therefore choose small constants $\eps, \delta > 0$ such that $s_1(t_0-\eps) = \delta$ and run the algorithm until stage $(t_0 - \eps) n$.  At this point $(t_0-\eps) n$ literals are already implied, proving the lower bound of~\eqref{eq:implied}.  

To prove the upper bound of~\eqref{eq:implied}, recall that by assumption $t_0$ is a simple root of~\eqref{eq:1}, i.e., the second derivative of the left-hand size of~\eqref{eq:1} with respect to $t$ is nonzero at $t_0$.  It is easy to see that this is equivalent to $\ds_1/\dt < 0$ at $t_0$.  Since $\ds_1/\dt$ is analytic, there is a constant $c > 0$ such that $\ds_1/\dt < 0$ for all $t_0 - c \le t \le t_0 + c$.  Set $\eps < c$; the number of literals implied during these stages is bounded above by a subcritical branching process whose initial population is \whp\ $\delta n + o(n)$, and by standard arguments we can bound its total progeny to be $\eps' n$ for any $\eps' > 0$ by taking $\delta$ small enough.
\remove{
To obtain the second part of inequality~\ref{eq:implied} we bound the 
number of literals labeled as implied after stage $T_\delta$. 
\begin{claim}
For every $\delta^{\prime}>0$ one can choose a constant $\delta >0$ such that 
\whp\ the number of literals labeled implied {\em after} stage $T_\delta$ 
is at most $\delta' \cdot n$. 
\end{claim}
We outline the proof of the previous claim (details are left for the full 
version). Note that the function $s_1$ is monotonically decreasing in the 
neighbourhood of $t_0$, since $s_1$ is analytical, $s_1(0)=d_1>0$ and $t_0$ is the smallest root of $s_1$. Therefore $s_1^{\prime}(t_0)\le 0$ and, since $t_0$ is not critical, $\lambda := s_1^{\prime}(t_0)< 0$. We infer that there exists a neighbourhood $(t_1,t_2)$ of $t_0$ such that 
$\forall t\in (t_1,t_2), s_1^{\prime}(t)<\lambda/2$.   

Then the number of literals implied at stages between $t_1 n$ and $t_2 n$ is  then dominated by a branching process with $\delta n$ initial individuals and expected birth rate $1+\lambda/2 < 1$. Choosing $\delta$ small enough we can make the number of individuals in the dominating branching process that are 
born before extinction is \whp\ at most $\delta' \cdot n$. 
}
\end{proof} 

It is easy to see that the backbone of implied positive literals is a uniformly random subset of $\{x_1,\ldots,x_n\}$ of size $N_{\boldd,n}$.  Since $x_1$ is guaranteed to not be among the $d_1 n$ initially positive literals, the probability that $x_1$ is not in this backbone is
\[ 
\frac{n-N_{\boldd,n}}{(1-d_1)\cdot n}.   
\]
By completeness of positive unit resolution for Horn satisfiability, this is 
precisely the probability that the $\phi$ is satisfiable. 
Applying Lemma~\ref{lem:gen} and taking $\epsilon \goesto 0$ proves equation~\eqref{eq:prsat} and completes the proof of Theorem~\ref{thm:main}.

We make several observations.  First, if we set $k=2$ and take the limit $d_1 \to 0$, 
Theorem~\ref{lem:gen} recovers Karp's result {\rm~\cite{Kar90}} on the size 
of the reachable component of a random directed graph with mean out-degree 
$d=d_2$, the root of $\ln (1-t) + d t = 0$.  

Secondly, as we will see below, the condition that $t_0$ is simple is essential.  Indeed, for the 1-3-Horn-SAT model studied in~\cite{DV05}, the curve $\Gamma$ of discontinuities, where the probability of satisfiability drops in the ``waterfall'' of Figure~\ref{fig:1-3}, consists exactly of those $(d_1,d_3)$ where $t_0$ is a double root, which implies $\ds_1/\dt = 0$ at $t_0$.

Finally, we note that Theorem~\ref{lem:gen} is very similar to Darling and Norris's work on 
identifiability in random undirected hypergraphs~\cite{DN05}, where the number of hyperedges of length $j$ is Poisson-distributed with mean $\beta_j$.  Their result reads
\[ \ln (1-t) + \sum_{j=1}^{k} j \beta_j t^{j-1} = 0 \enspace . \]
We can make this match~\eqref{eq:1} as follows.  First, since each hyperedge of length $j$ corresponds to $j$ Horn clauses, we set $d_j = j \beta_j$ for all $j \ge 2$.  Then, since edges are chosen with replacement in their model, the expected number of distinct clauses of length $1$ (i.e., positive literals) is $d_1 n$ where $d_1 = 1-\e^{-\beta_1}$.

\section{Application to $H_{n,\boldd}^2$}

\begin{figure}[tbh]
\begin{center} 
\begin{minipage}[t]{2.2 in}
\begin{center}
\psfig{figure=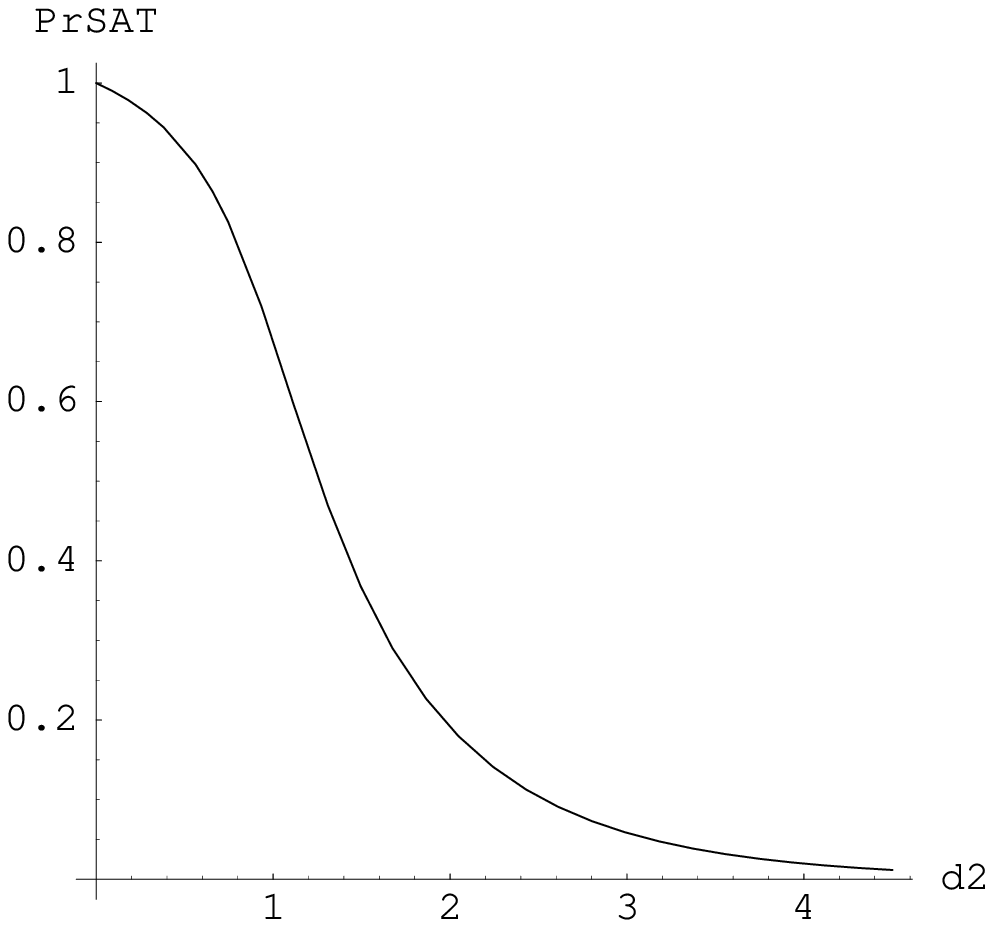,height=50truemm,width=50truemm}
\end{center}
\end{minipage}
\begin{minipage}[t]{2.3 in}
\begin{flushright}
\hspace{0.5in}\psfig{figure=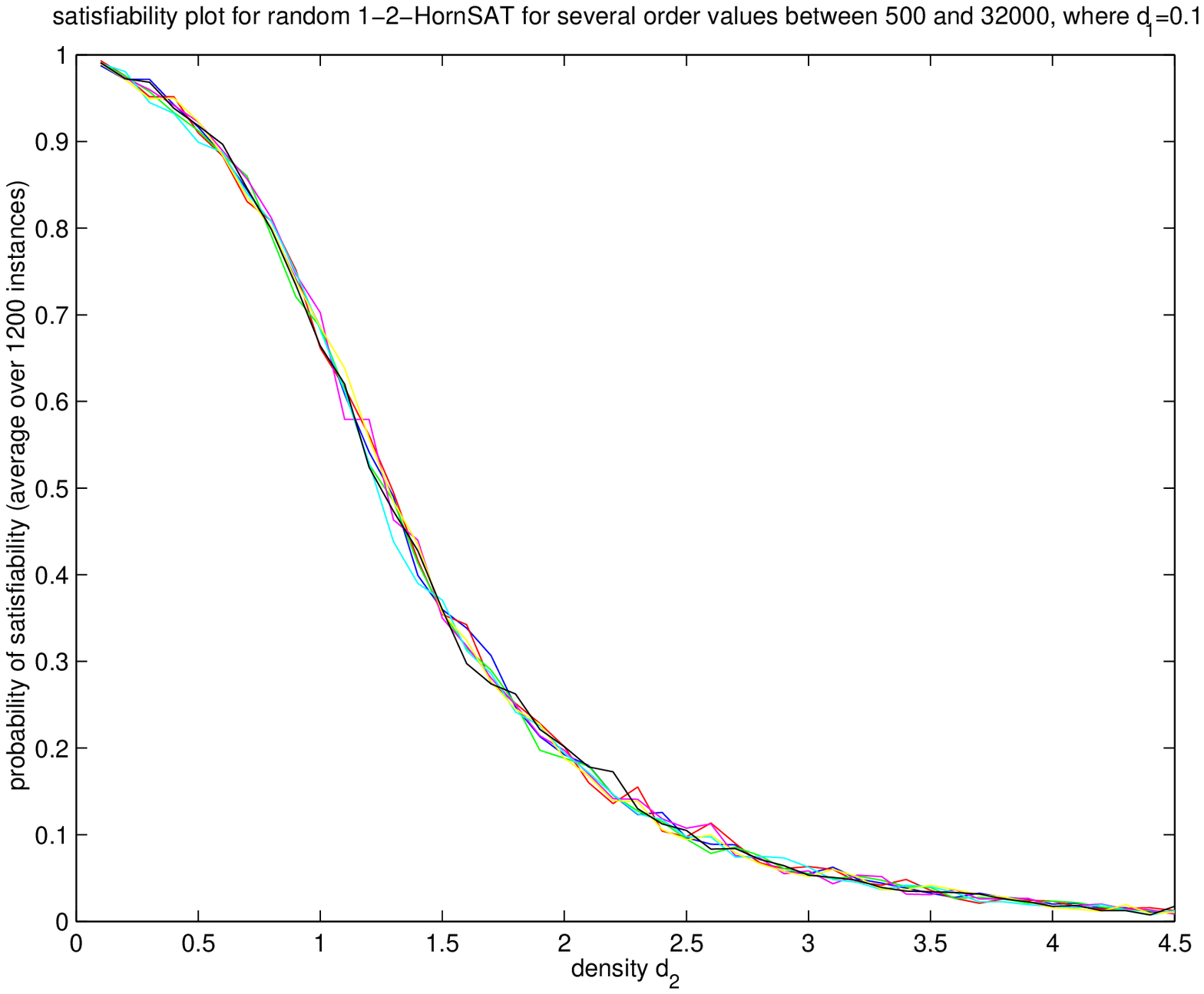,height=50truemm,width=45truemm} 
\end{flushright}
\end{minipage}
\caption{The probability of satisfiability for 1-2-Horn formulae as a function 
of $d_2$, where $d_1 = 0.1$.  Left, our analytic results; right, the experimental data of~\cite{DV05}.}  
\label{1-2-0.1}
\end{center} 
\end{figure}

For $H_{n,\boldd}^2$, Equation~\eqref{eq:1} can rewritten 
as $1-t=(1-d_1)\cdot e^{-d_2\cdot t}$. Denoting $y=d_2(t-1)$,
this implies 
\[
y\cdot e^{y} = d_2(t-1)\cdot e^{d_2\cdot (t-1)}= -d_2(1-d_1)\cdot e^{-d_2\cdot t}\cdot e^{d_2\cdot (t-1)}=-(1-d_1)d_2 e^{-d_2}.
\]
Solving this yields  
\[ t_0 = 1 + \frac{1}{d_2} W\bigl( -(1-d_1) d_2 \e^{-d_2} \bigr) \]
and substituting this into~\eqref{eq:prsat} 
\remove{
\begin{equation}
\label{eq:prop2}
\Phi(d_1,d_2) 
= \lim_{n \to \infty} \Pr[H^2_{n,\boldd} \mbox{ is satisfiable}] = 
-\frac{W\bigl( -(1-d_1) d_2 \e^{-d_2} \bigr)}{(1-d_1)d_2},
\end{equation}
}
proves Equation~\eqref{eq:prop2} and Proposition~\ref{prop:2}.  In Figure~\ref{1-2-0.1}
we plot the probability of satisfiability $\Phi(d_1,d_2)$ as a function of $d_2$ for $d_1=0.1$ and compare it with the experimental results of~\cite{DV05}; the agreement is excellent.

\section{A continuous-discontinuous phase transition for 
$H_{n,d_1,0,d_3}^3$}

\begin{figure}[tbh]
\begin{center} 
\begin{minipage}[t]{3 in}
\begin{center}
\psfig{figure=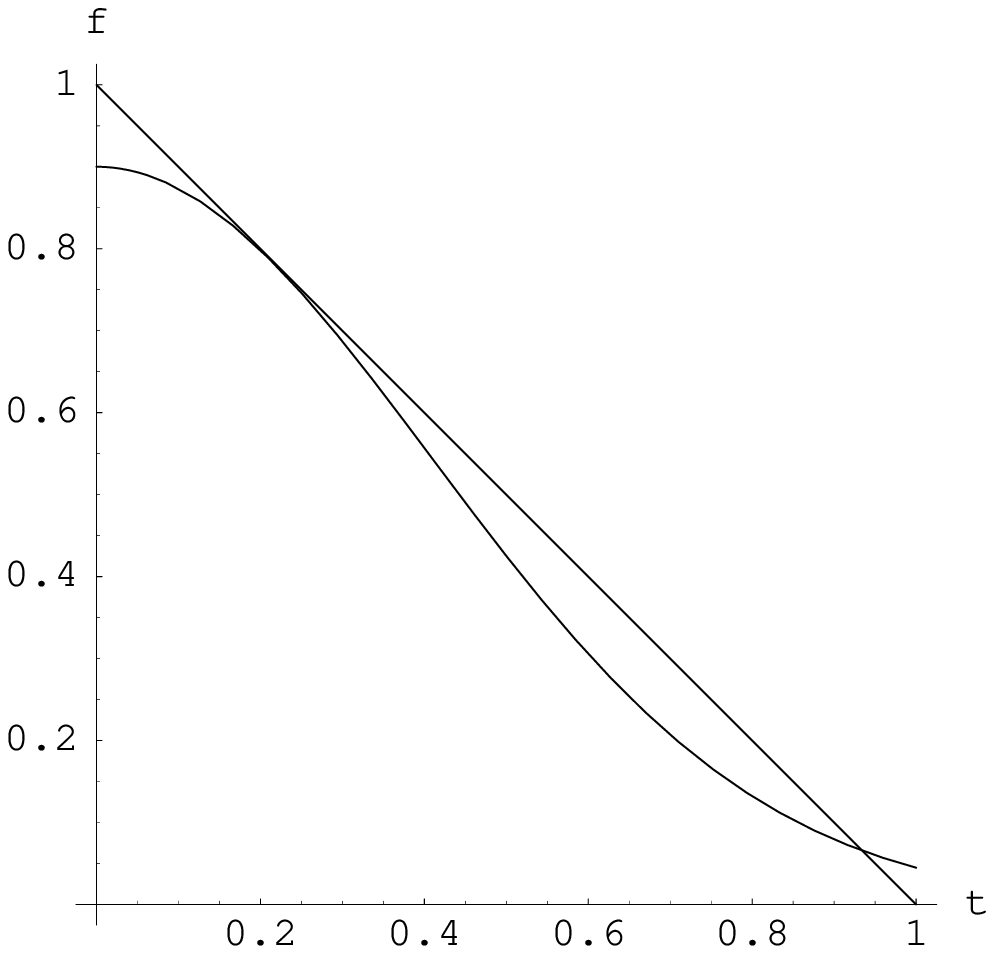,height= 60 truemm,width= 60 truemm}
\end{center}
\end{minipage}
\begin{minipage}[t]{3 in}
\begin{flushright}
\hspace{0.25in}\psfig{figure=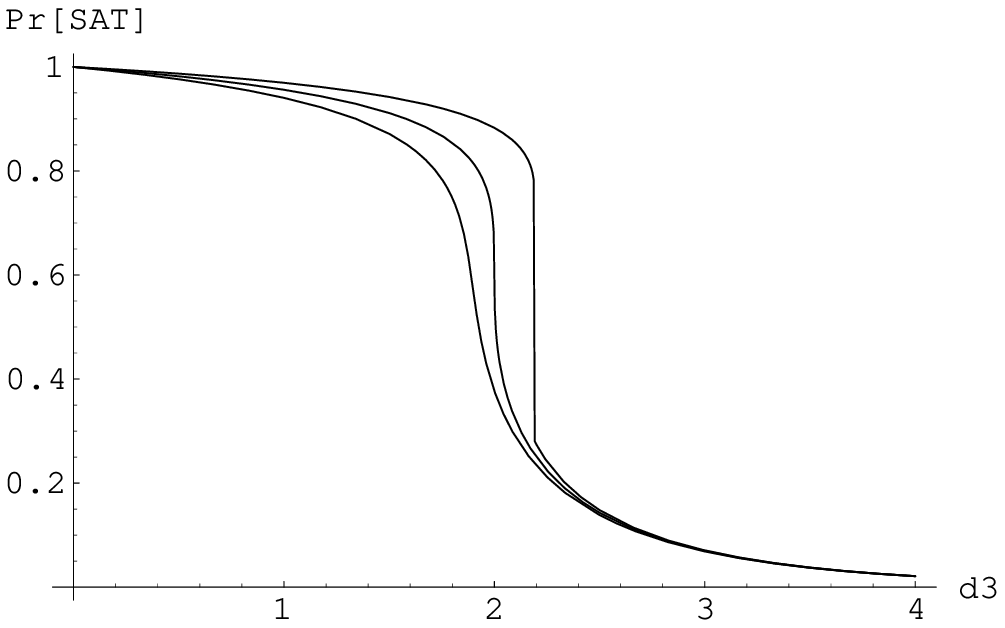,height=60 truemm, width= 65 truemm} 
\end{flushright}
\end{minipage}
\caption{Left, the function $f(t)$ of~\eqref{eq:f} when $d_1 = 0.1$ and $d_3 = 3$. 
Right, the probability of satisfiability $\Phi(d_1,d_3)$ with $d_1$ equal to  
$0.15$ (continuous), $0.1756$ (critical), and $0.2$ (discontinuous).}  
\label{first-f}
\end{center}
\end{figure}

For the random model $H_{n,d_1,0,d_3}^3$ studied in~\cite{DV05}, 
an analytic solution analogous to~\eqref{eq:prop2} does not seem to exist.  
Let us, however, rewrite~\eqref{eq:1} as
\begin{equation}
\label{eq:f}
 1-t = f(t) := (1-d_1) \e^{-d_3 t^2} \enspace . 
\end{equation}
We claim that for some values of $d_1$ and $d_3$ there is a phase transition in the roots of~\eqref{eq:f}.  For instance, consider the plot of $f(t)$ shown in Figure~\ref{first-f} for $d_1 = 0.1$ and $d_3 = 3$. Here $f(t)$ is tangent to $1-t$, so there is a bifurcation as we vary either parameter; for $d_3 = 2.9$, for instance, $f(t)$ crosses $1-t$ three times and there is a root of~\eqref{eq:f} at $t=0.185$, but for $d_3 = 3.1$ the unique root is at $t=0.943$.  This causes the probability of satisfiability to jump discontinuously but from $0.905$ to $0.064$.  The set of pairs $(d_1,d_3)$ for which $f(t)$ is tangent to $1-t$ is exactly the curve $\Gamma$ on which the smallest positive root $t_0$ of~\eqref{eq:1} or~\eqref{eq:f} is a double root, giving the ``waterfall'' of Figure~\ref{fig:1-3}.

To find where this transition takes place, we set $f' = -1$, yielding $f(t) = 1/(2 d_3 t)$. Setting this to $1-t$ and solving for $t$ gives
\begin{equation}
\label{eq:d1a}
 d_1 = 1 - \frac{\e^{d_3 t^2}}{2 d_3 t} 
\end{equation}
where
\begin{equation}
\label{eq:t}
 t = \frac{1}{2} \left( 1 - \sqrt{1-\frac{2}{d_3}} \right) 
\end{equation}
Substituting~\eqref{eq:t} in~\eqref{eq:d1a} and simplifying gives~\eqref{eq:gamma}, 
\remove{
\begin{equation}
\label{eq:d1}
 d_1 = 1 - \frac{\exp\left( \frac{1}{4} \left( \sqrt{d_3} - \sqrt{d_3-2} \right)^2 \right)}
{d_3 - \sqrt{d_3 (d_3-2)}},
\end{equation}
}
proving Proposition~\ref{prop:3}.

The fact that $d_1$ is only real for $d_3 \ge 2$ explains why $\Gamma$ 
ends at $d_3 = 2$.  At this extreme case we have 
\[ d_1 = 1-\frac{\sqrt{\e}}{2} \approx 0.1756 \;\mbox{ and }\;
\frac{\partial d_1}{\partial d_3} = - \frac{\sqrt{\e}}{8} \enspace . \]
What happens for $d_3 < 2$?  In this case, there are no real $t$ for which 
$f'(t) = -1$, so the kind of tangency displayed in 
Figure~\ref{first-f} cannot happen.  
In that case,~\eqref{eq:f} (and equivalently~\eqref{eq:1}) has a unique 
solution $t$, which varies continuously with $d_1$ and $d_3$, and therefore 
the probability of satisfiability $\Phi(d_1,d_3)$ is continuous as well. 
To illustrate this, in the right part of Figure~\ref{first-f} we plot
$\Phi(d_1,d_3)$ as a function of $d_3$ with three values of $d_1$.  
For $d_1=0.15$, $\Phi$ is continuous; for $d_1 = 0.2$, it is discontinuous; 
and for $d_1 = 0.1756...$, the critical value at the second-order transition, 
it is continuous but has an infinite derivative at $d_3=2$.

\newcommand{\etalchar}[1]{$^{#1}$}


{\small
\begin{thebibliography}{MZK{\etalchar{+}}99}

\bibitem{achreview}
D.~Achlioptas.
\newblock Lower Bounds for Random 3-{SAT} via Differential Equations.
\newblock {\em Theoretical Computer Science} 265 (1-2), 159--185, 2001.

\bibitem{ACIM01}
D.~Achlioptas, A.~Chtcherba, G.~Istrate, and C.~Moore.
\newblock The phase transition in 1-in-$k$ {SAT} and {NAE} 3-{SAT}.
\newblock In {\em Proc.~12th ACM-SIAM Symp. on Discrete Algorithms}, 721--722, 2001.

\bibitem{AKKK01}
D. Achlioptas, L.M. Kirousis, E. Kranakis, and D. Krizanc.
\newblock Rigorous results for random $(2+p)$-SAT.
\newblock {\em Theor. Comput. Sci.} 265(1-2): 109-129 (2001).

\bibitem{BCD01}
T.~Bench-Capon and P.~Dunne.
\newblock A sharp threshold for the phase transition of a restricted
  satisfiability problem for {H}orn clauses.
\newblock {\em Journal of Logic and Algebraic Programming}, 47(1):1--14, 2001.

\bibitem{BDFN92}
J. J. Binney, N. J. Dowrick, A. J. Fisher and M. E. J. Newman.
{\em The Theory of Critical Phenomena.} Oxford University Press (1992).

\bibitem{Bol85}
B.~Bollob\'{a}s.
\newblock {\em Random Graphs}.
\newblock Academic Press, 1985.

\bibitem{CDMM03}
S.~Cocco, O.~Dubois, J.~Mandler, and R.~Monasson.
\newblock Rigorous decimation-based construction of ground pure states
for spin glass models on random lattices.
\newblock {\em Phys. Rev. Lett.}, 90, 2003.

\bibitem{CR92}
V.~Chv\'{a}tal and B.~Reed.
\newblock Mick gets some (the odds are on his side).
\newblock In {\em Proc.~33rd IEEE Symp. on Foundations of Computer
Science}, IEEE Comput. Soc. Press, 620--627, 1992.


\bibitem{CA96}
J.~M. Crawford and L.~D. Auton.
\newblock Experimental results on the crossover point in random 3-{SAT}.
\newblock {\em Artificial Intelligence}, 81(1-2):31--57, 1996.

\bibitem{DG84}
W.~F. Dowling and J.~H. Gallier.
\newblock Linear-time algorithms for testing the satisfiability of
propositional {H}orn formulae.
\newblock {\em Logic Programming. (USA) ISSN: 0743-1066}, 1(3):267--284,
1984.

\bibitem{DBM00}
O.~Dubois, Y.~Boufkhad, and J.~Mandler.
\newblock Typical random 3-{SAT} formulae and the satisfiability theshold.
\newblock In {\em Proc.~11th ACM-SIAM Symp. on Discrete Algorithms}, 126--127, 2000.

\bibitem{DM02}
O.~Dubois and J.~Mandler.
\newblock The 3-{XORSAT} threshold.
\newblock In {\em Proc.~43rd IEEE Symp. on Foundations of Computer Science}, 769--778, 2002.

\bibitem{DN05}
R.~Darling and J.R.~Norris.
\newblock Structure of large random hypergraphs.
\newblock {\em Annals of Applied Probability}, 15(1A), 2005.

\bibitem{DV05}
D.~Demopoulos and M.~Vardi.
\newblock The phase transition in random 1-3 {H}ornsat problems.
\newblock In A.~Percus, G.~Istrate, and C.~Moore, editors, {\em Computational
  Complexity and Statistical Physics}, Santa Fe Institute Lectures in the
  Sciences of Complexity. Oxford University Press, 2005.
\newblock Available at {\tt http://www.cs.rice.edu/$^\sim$vardi/papers/}.

\bibitem{ER60}
P.~Erd\"{o}s and A.~R\'{e}nyi.
\newblock On the evolution of random graphs.
\newblock {\em Publications of the Mathematical Institute of the Hungarian Academy of Science}, 5:17--61, 1960.

\bibitem{Fri99}
E.~Friedgut.
\newblock Necessary and sufficient conditions for sharp threshold of graph
  properties and the $k$-{SAT} problem.
\newblock {\em J. Amer. Math. Soc.}, 12:1917--1054, 1999.

\bibitem{Goe96}
A.~Goerdt.
\newblock A threshold for unsatisfiability.
\newblock {\em J. Comput. System Sci.}, 53(3):469--486, 1996.

\bibitem{HajiaSorkin}
M.~Hajiaghayi and G.B.~Sorkin.
\newblock The satisfiability threshold for random 3-SAT is at least 3.52.
\newblock IBM Technical Report , 2003.

\bibitem{HW94}
T.~Hogg and C.~P. Williams.
\newblock The hardest constraint problems: A double phase transition.
\newblock {\em Artificial Intelligence}, 69(1-2):359--377, 1994.

\bibitem{Ist02}
G.~Istrate.
\newblock The phase transition in random {H}orn satisfiability and its
  algorithmic implications.
\newblock {\em Random Structures and Algorithms}, 4:483--506, 2002.

\bibitem{Ist04}
G.~Istrate.
\newblock On the satisfiability of random $k$-{H}orn formulae.
\newblock In P.~Winkler and J.~Nesetril, editors, {\em Graphs, Morphisms and
  Statistical Physics}, volume~64 of {\em AMS-DIMACS Series in Discrete
  Mathematics and Theoretical Computer Science}, 113--136, 2004.

\bibitem{Kar90}
R.~Karp.
\newblock The transitive closure of a random digraph.
\newblock {\em Random Structures and Algorithms}, 1:73--93, 1990.

\bibitem{KKL}
A.~Kaporis, L.~M. Kirousis, and E.~Lalas.
\newblock Selecting complementary pairs of literals.
\newblock In {\em Proceedings of LICS'03 Workshop
on Typical Case Complexity and Phase Transitions}, June 2003.

\bibitem{Mak87}
J.~A. Makowsky.
\newblock Why {H}orn formulae matter in {C}omputer {S}cience: Initial
  structures and generic examples.
\newblock {\em JCSS}, 34(2-3):266--292, 1987.

\bibitem{MZ02}
M.~M\'{e}zard and R.~Zecchina.
\newblock Random k-satisfiability problem: from an analytic solution to an efficient algorithm.
\newblock {\em Phys. Rev. E}, 66:056126, 2002.

\bibitem{MPZ02}
M.~M\'{e}zard, G.~Parisi, and R.~Zecchina.
\newblock Analytic and algorithmic solution of random satisfiability problems.
\newblock {\em Science}, 297:812--815, 2002.

\bibitem{MZKST99}
R.~Monasson, R.~Zecchina, S.~Kirkpatrick, B.~Selman, and L.~Troyansky.
\newblock $2+p$-{SAT}: Relation of typical-case complexity to the nature of the
  phase transition.
\newblock {\em Random Structures and Algorithms}, 15(3--4):414--435, 1999.

\bibitem{SML96}
B.~Selman, D.~G. Mitchell, and H.~J. Levesque.
\newblock Generating hard satisfiability problems.
\newblock {\em Artificial Intelligence}, 81(1-2):17--29, 1996.

\bibitem{SK96}
B.~Selman and S.~Kirkpatrick.
\newblock Critical behavior in the computational cost of satisfiability testing.
\newblock {\em Artificial Intelligence}, 81(1-2):273--295, 1996.

\bibitem{VW86a}
M.Y. Vardi and P.~Wolper.
\newblock Automata-Theoretic Techniques for Modal Logics of Programs
\newblock {\em J. Computer and System Science} 32:2(1986), 181--221, 1986.

\bibitem{VW86b}
M.Y. Vardi and P.~Wolper.
\newblock An automata-theoretic approach to automatic program 
verification (preliminary report).
\newblock In {\em Proc.~1st IEEE Symp. on Logic in Computer Science}, 332--344, 1986.

\bibitem{Wor95}
N.~Wormald.
\newblock Differential equations for random processes and random graphs.
\newblock {\em Annals of Applied Probability}, 5(4):1217--1235, 1995.


\end{thebibliography}
}
\end{document}